\documentclass[10pt]{article}

\date{} 
 
\usepackage{hyperref}
\usepackage{amssymb, amsmath}
\def\sideremark#1{\ifvmode\leavevmode\fi\vadjust{\vbox to0pt{\vss
 \hbox to 0pt{\hskip\hsize\hskip1em
\vbox{\hsize2cm\small\raggedright\pretolerance10000eb
 \noindent #1\hfill}\hss}\vbox to8pt{\vfil}\vss}}}

\usepackage{amsmath}
\usepackage{latexsym}
\usepackage{amssymb}

\usepackage{amsthm}

\font\tengoth=eufm10 at 10pt
\font\sevengoth=eufm7 at 6pt
\newfam\gothfam
\textfont\gothfam=\tengoth
\scriptfont\gothfam=\sevengoth

\newcommand{\mlabel}[1]{\label{#1}}

\renewcommand{\:}{\colon}

\newcommand{\cE}{\mathcal{E}}

\newcommand{\cH}{\mathcal{H}}

\newcommand{\cL}{\mathcal{L}}

\newcommand{\cN}{\mathcal{N}}

\newcommand{\subeq}{\subseteq}

\newcommand{\N}{{\mathbb N}}
\newcommand{\Z}{{\mathbb Z}}
\newcommand{\R}{{\mathbb R}}
\newcommand{\C}{{\mathbb C}}

\newcommand{\Q}{{\mathbb Q}}

\renewcommand{\hat}{\widehat}

\renewcommand{\tilde}{\widetilde}

\newcommand{\id}{\mathop{{\rm id}}\nolimits}

\newcommand{\vphi}{\varphi}
\renewcommand{\phi}{\varphi}

\newcommand{\Rarrow}{\Rightarrow}
 
\newcommand{\oline}{\overline}

\newcommand{\Larrow}{\Leftarrow}
\newcommand{\la}{\langle}
\newcommand{\ra}{\rangle}

\newcommand{\res}{\vert}

\newcommand{\ssssarr}{\hbox to 15pt{\rightarrowfill}}
\newcommand{\sssarr}{\hbox to 20pt{\rightarrowfill}}
\newcommand{\ssarr}{\hbox to 30pt{\rightarrowfill}}
\newcommand{\sarr}{\hbox to 40pt{\rightarrowfill}}
\newcommand{\arr}{\hbox to 60pt{\rightarrowfill}}
\newcommand{\larr}{\hbox to 60pt{\leftarrowfill}}
\newcommand{\Arr}{\hbox to 80pt{\rightarrowfill}}

\def\theoremname{Theorem}
\def\propositionname{Proposition}
\def\corollaryname{Corollary}
\def\lemmaname{Lemma}
\def\remarkname{Remark}
\def\conjecturename{Conjecture} 

\def\definitionname{Definition}
\def\exercisename{Exercise}
\def\examplename{Example}
\def\examplesname{Examples}
\def\problemname{Problem}
\def\problemsname{Problems}

\def\satzname{Satz} 
\def\koroname{Korollar}
\def\folgname{Folgerung}
\def\bemerkname{Bemerkung}
\def\aufgname{Aufgabe}

\def\beisname{Beispiel}
\def\beissname{Beispiele}
\def\bewname{Beweis}

\def\@thmcounter#1{\noexpand\arabic{#1}}
\def\@thmcountersep{}
\def\@begintheorem#1#2{\it \trivlist \item[\hskip 
\labelsep{\bf #1\ #2.\quad}]}
\def\@opargbegintheorem#1#2#3{\it \trivlist
      \item[\hskip \labelsep{\bf #1\ #2.\quad{\rm #3}}]}
\makeatother
\newtheorem{theor}{\theoremname}[section]
\newtheorem{propo}[theor]{\propositionname}
\newtheorem{coro}[theor]{\corollaryname}
\newtheorem{lemm}[theor]{\lemmaname}

\newenvironment{thm}{\begin{theor}\it}{\end{theor}}
\newenvironment{theorem}{\begin{theor}\it}{\end{theor}}

\newenvironment{prop}{\begin{propo}\it}{\end{propo}}

\newenvironment{cor}{\begin{coro}\it}{\end{coro}}

\newenvironment{lem}{\begin{lemm}\it}{\end{lemm}}

\newtheorem{rema}[theor]{\remarkname}

\newenvironment{remark}{\begin{rema}\rm}{\end{rema}}
\newenvironment{rem}{\begin{rema}\rm}{\end{rema}}

\newtheorem{stepnow}[theor]{}

\newtheorem{defin}[theor]{\definitionname} %% write

\newenvironment{definition}{\begin{defin}\rm}{\end{defin}}
\newenvironment{defn}{\begin{defin}\rm}{\end{defin}}

\newtheorem{exerc}{\exercisename}[section]

\newtheorem{exa}[theor]{\examplename}

\newenvironment{example}{\begin{exa}\rm}{\end{exa}}
\newenvironment{ex}{\begin{exa}\rm}{\end{exa}}

\newtheorem{exas}[theor]{\examplesname}

\newenvironment{exs}{\begin{exas}\rm}{\end{exas}}

\newtheorem{conj}[theor]{\conjecturename}

\newtheorem{pro}[theor]{\problemname}

\newtheorem{prs}[theor]{\problemsname}

\newtheorem{aufg}{\aufgname}[section]

\newenvironment{prf}{\begin{proof}}{\end{proof}}
 
\newcommand{\pmat}[1]{\begin{pmatrix} #1 \end{pmatrix}}

{\hfill\qed\end{trivlist}}

\newenvironment{beweis*}{\begin{trivlist}\item[\hskip%
\labelsep{\bf\bewname.\quad}]}% 
{\end{trivlist}}

\newtheorem{satzn}[theor]{\satzname}

\newtheorem{koro}[theor]{\koroname}

\newtheorem{folg}[theor]{\folgname}

\newtheorem{bem}[theor]{\bemerkname}

\newtheorem{aufgn}[theor]{\aufgname}

\newtheorem{beis}[theor]{\beisname}

\newtheorem{beiss}[theor]{\beissname}

\begin{document}
%%%%%%%%%%%%%%%%%%%%%%%%%%

\title{Reflection positivity on real intervals}
% \\ {\tt interval.tex}} 
\author{P. Jorgensen, K.-H. Neeb, G. \'Olafsson} 

\maketitle

\begin{abstract}
We study functions $f \: (a,b) \to \R$ on open intervals in 
$\R$ with respect to various kinds of positive and negative 
definiteness conditions. We say that 
$f$ is positive definite if the kernel 
$f\big(\frac{x + y}{2}\big)$ is positive definite. We call 
$f$ negative definite if, for every $h > 0$, the function 
$e^{-hf}$ is positive definite. 
Our first main result is a L\'evy--Khintchine formula 
(an integral representation) for negative definite functions 
on arbitrary intervals. For $(a,b) = (0,\infty)$ it generalizes 
classical results by Bernstein and Horn. 

On a symmetric interval $(-a,a)$, we call $f$ reflection positive 
if it is positive definite and, in addition, the kernel
$f\big(\frac{x - y}{2}\big)$ is positive definite. We likewise define 
reflection negative functions and obtain a 
L\'evy--Khintchine formula for reflection negative functions 
on all of $\R$. 
Finally, we obtain a characterization of germs of reflection 
negative functions on $0$-neighborhoods in $\R$. 
\end{abstract}

\section{Introduction}

Positive definiteness conditions for functions 
$f \: (a,b) \to \C$ on intervals $(a,b) \subeq \R$ are a classical 
subject of analysis, operator theory, stochastics and 
harmonic analysis (\cite{Wi46},  \cite{SzNBK10}, \cite{SSV10},\cite{BCR84}). 
Basically, there are two types of positive definiteness conditions. 
The first one comes from the additive group $(\R,+)$, for which a 
function $f \: \R \to \C$ is positive definite if and only if the kernel 
$(f(x-y))_{x,y \in \R}$ is positive definite. This condition makes also sense 
on symmetric intervals of the form $(-a/2,a/2)$ if $f$ is defined on 
$(-a,a)$ (see \cite{JN15, JPT15} for recent progress in the local theory). 
Bochner's Theorem asserts that a continuous positive definite function 
on $(\R,+)$ is the Fourier transform 
$f(t) = \int_\R e^{-it\lambda}\, d\mu(\lambda)$ of a bounded positive measure 
$\mu$ on~$\R$. 

The second type makes sense for functions 
$f \: (a,b) \to\C$ on any interval and requires that the kernel 
$\big(f\big(\frac{x+y}{2}\big)\big)_{a < x,y < b}$ is positive definite.  
By Widder's Theorem, this is equivalent to $f$ being a Laplace transform 
$f(t) = \int_\R e^{-t\lambda}\, d\mu(\lambda)$ of a positive measure $\mu$ on~$\R$.
For $(a,b) = (0,\infty)$ this is precisely the condition of 
positive definiteness on the $*$-semigroup $(0,\infty)$ with the trivial 
involution $t^* =t$ for $t > 0$. 

We call a function $f \:  (-a,a) \to \C$ on a symmetric interval 
{\it reflection positive} if it satisfies two positive definiteness 
conditions at the same time, namely if both kernels 
\[ \Big(f\Big(\frac{x-y}{2}\Big)\Big)_{a < x,y < b} \quad \mbox{ and } \quad 
\Big(f\Big(\frac{x+y}{2}\Big)\Big)_{a < x,y < b} \] 
are positive definite. This notion is motivated by our work 
on reflection positive representations of the additive 
group $(\R,+)$, where such functions 
occur for $a = \infty$ \cite{NO14} and for 
periodic functions \cite{NO15, NO16}. 

By Schoenberg's Theorem, any notion of positive definiteness determines a corresponding 
notion of negative definiteness: a function $\psi$ is called {\it 
negative definite} if, for each $h > 0$, the function 
$e^{-h\psi}$ is positive definite. Accordingly, we define 
reflection negative functions. Originally, the present article 
was motivated by the reflection negative functions on 
$\R$ arising in \cite{JNO16} in the context of affine reflection positive 
actions of $\R$ on affine Hilbert spaces corresponding to Gaussian 
processes with stationary square increments. 

While the notion of a positive definite function and related extension questions
make sense in a very general context, we shall restrict attention here to the case
when the domain is an interval. Our main focus here is the framework of
reflection positivity as it occurs in Quantum Field Theory.
As we point out below, reflection positivity is of great interest in
various non-commutative settings. However, even for the seemingly
``easy'' case of a function on an interval, the study of  reflection positivity
is non-trivial and of independent interest. Moreover,
many questions in more general settings are more transparent in the simpler case of intervals.
While we shall not presently follow up on all connections to neighbouring
fields, we shall hint  here at
intriguing connections between Bernstein functions and operator monotonocity
(see e.g., \cite[Ch.~11]{SSV10}), connections with
the theory of non-commutative Pick-interpolation theory
\cite{Po08, AM15}  and connections to
operator monotone functions in several commuting variables \cite{AMY12}
and \cite{Pa14} in the case of non-commuting variables.

The main new results of the present paper are characterizations, 
resp., integral represetations of the following classes of functions: 
\begin{itemize}
\item Negative definite functions on arbitrary intervals 
(the kernel $\psi\big(\frac{x+y}{2}\big)$ is negative definite): 
Fix $t_0 \in (a,b)$.  Then $\psi \: (a,b) \to \R$ is negative definite 
if and only if there exists a 
positive measure $\mu$ on $\R$ and $c,d \in \R$ such that 
\[ \psi(t) = c + d(t-t_0) + \int_{-\infty}^\infty e_\lambda(t)e^{-\lambda t_0}
\, d\mu(\lambda), \ \  \mbox{ where } \ \ 
e_\lambda(t) := \begin{cases}
\frac{1 - \lambda (t-t_0) - e^{-\lambda (t-t_0)}}{\lambda^2} 
& \text{ for } \lambda\not=0 \\ 
-(t-t_0)^2/2  & \text{ for } \lambda=0.  
\end{cases}\] 
This result (Theorem~\ref{thm:L-K-interval}) generalizes classical results by 
Horn for $(a,b) = (0,\infty)$ \cite{Ho67} and  
Bernstein for $(a,b) = (0,\infty)$ and $\psi \geq 0$. 
\item  Increasing negative definite functions on 
$(0,\infty)$ (Theorem~\ref{thm:L-K-incr}) (this is closely related to Bernstein's 
Theorem~\ref{thm:lk-bernstein2}
 characterizing all non-negative negative definite functions):
A function $\psi \: (0, \infty) \to \R$ is negative definite and increasing 
if and only if there exists a positive measure $\mu$ on $\R$ and $c\in \R$ such that 
\begin{equation}
%  \label{eq:leki-inc}
\psi(t) = c + \int_{[0,\infty)} f_\lambda(t) \, d\mu(\lambda),
\quad \mbox{ where } \quad 
f_\lambda(t) := 
\begin{cases}
\frac{e^{-\lambda} - e^{-\lambda t}}{\lambda} 
& \text{ for } \lambda\not=0 \\ 
t-1  & \text{ for } \lambda=0.  
\end{cases}
\end{equation}
Then $c = \psi(1)$ and $\mu$ is uniquely determined by $\psi' = \cL(\mu)$. 
A positive measure $\mu$ 
occurs if and only if its Laplace transform $\cL(\mu)$ is finite on $(0,\infty)$. 
\item Functions which are reflection positive on some $0$-neighborhood in $\R$ 
(Theorem~\ref{thm:5.9}): 
Let $\mu$ be a finite positive measure on $\R$ for which 
$\vphi(t) := \cL(\mu)(|t|)$ exists for $|t| \leq a$. 
\begin{itemize}
\item[\rm(a)] If $\cL(\mu)'(a-) \leq 0$, then $\vphi$ is reflection positive on $[-a,a]$ 
and extends to a positive definite function on $\R$. 
\item[\rm(b)] If $\vphi$ is reflection positive on $[-a,a]$ and non-constant, 
then there exists an element $b \in (0,a]$ with $\cL(\mu)'(b-) < 0$. 
\end{itemize}
This relates naturally to the 
integral representations of reflection positive function 
on $\R$ in \cite{NO14} and of $\beta$-periodic reflection positive 
functions on $[-\beta,\beta]$ in \cite{NO15}.\\ 
\item Reflection negative functions on $\R$ 
(Theorem~\ref{thm:lk-bernstein}): 
A symmetric continuous function \break $\psi \colon \R \to \R$ is reflection 
negative if and only 
if $\psi\res_{(0,\infty)}$ is a Bernstein function. In particular, this is equivalent to 
the existence of $a,b \geq 0$ and a positive measure $\mu$ on $(0,\infty)$ 
with $\int_0^\infty (1 \wedge \lambda)\, d\mu(\lambda) < \infty$ 
such that 
\[ \psi(t) = a + b|t| + \int_0^\infty (1 - e^{-\lambda |t|})\, d\mu(\lambda). \] 
Here $a,b$ and $\mu$ are uniquely determined by~$\psi$.
\end{itemize}

This note is part of our long term project on \textit{reflection positivity}. 
This is a basic concept in constructive quantum 
field theory \cite{GJ81, JO00}, 
where it arises as a requirement on the euclidean side to establish a 
duality between euclidean and relativistic quantum field theories \cite{OS73}. 
The notion of a reflection negative function stems from 
\cite{JNO16}, where we study reflection positivity for affine actions 
of on a real Hilbert space. 

\tableofcontents 

\section{Reflection positive kernels}

We first recall the basic definitions concerning 
positive and negative definite kernels. 
As customary in physics, 
we follow the convention that the inner product of a complex
 Hilbert space is linear in the second argument.

\begin{defn}   \label{def:8.1.1} 
(a) Let $X$ be a set. A kernel $Q \colon X \times X \to \C$ is called 
{\it hermitian} if $Q(x,y) = \overline{Q(y,x)}$. 
A hermitian kernel $Q$ is called
 {\it positive definite} if 
$\sum_{j,k=1}^n c_j \overline{c_k} Q(x_j, x_k) \geq 0$ 
holds for $x_1, \ldots, x_n \in X, c_1, \ldots, c_n \in \C$. 
It is called {\it negative definite} if 
$\sum_{j,k=1}^n c_j \overline{c_k} Q(x_j, x_k) \leq 0$
holds for $x_1, \ldots, x_n \in X$ and $c_1, \ldots, c_n \in \C$ 
satisfying $\sum_{j = 1}^n c_j=0$ (\cite{BCR84}). 

(b) If $G$ is a group, then a function $\vphi \colon G \to \C$ is called {\it positive (negative) definite} if the kernel $(\vphi(gh^{-1}))_{g,h \in G}$ is positive (negative) definite. 
More generally, if $(S,*)$ is an involutive semigroup, then $\vphi \colon S \to \C$ is 
called {\it positive (negative) definite} if the kernel $(\vphi(st^*))_{s,t\in S}$ 
is positive (negative) definite. 
\end{defn} 

\begin{remark} We point out that a function $\psi :G \to \C$ is negative definite if and only if, for every $h>0$, the function
$e^{-h\psi}$ is positive definite (\cite[Thm. 3.2.2.]{BCR84}.
\end{remark}

\begin{remark} \label{rem:kerspace} 
Let $X$ be a set, $K \colon X \times X \to \C$ be a positive definite 
kernel and $\cH_K \subeq \C^X$ be the corresponding {\it reproducing kernel 
Hilbert space}. This is the unique Hilbert subspace of $\C^X$ on which all 
point evaluations $f \mapsto f(x)$ are continuous and given by 
\[ f(x) = \la K_x, f \ra \quad \mbox{ for } \quad 
K_y(x) = K(x,y)  := \la K_x, K_y \ra.\] 
\end{remark}

\begin{defn} \mlabel{def:1.2}  A {\it reflection positive Hilbert space} is a triple 
$(\cE,\cE_+,\theta)$, where 
$\cE$ is a Hilbert space, $\theta$ is a unitary involution and 
$\cE_+$ is a closed subspace which is {\it $\theta$-positive} in the sense that 
$\la \theta v, v\ra \geq 0$ for $v\in\cE_+$. 
For a reflection positive Hilbert space  $(\cE,\cE_+,\theta)$, 
let $\cN:=\{u \in\cE_+\: \la \theta u, u \ra =0\}$ and let
$\hat\cE$ be the completion of $\cE_+/\cN$ with respect to the inner product   
$\la \cdot, \cdot \ra_\theta$. We write 
 $q \: \cE_+ \to \hat\cE$ for the canonical map. 
\end{defn}

\begin{example} \mlabel{ex:1.3} 
(a) Suppose that $K \colon X \times X \to \C$ is a positive definite kernel, 
that $\tau \colon X \to X$ is an involution 
leaving $K$ invariant and that $X_+ \subeq X$ is a subset with the property that the 
kernel $K^\tau(x,y) := K(\tau x, y)= K(x, \tau y)$ is positive definite on $X_+$. 
Then the closed subspace $\cH_K^+ \subeq \cH_K$ generated 
by $(K_x)_{x \in X_+}$ is $\theta$-positive for 
$(\theta f)(x) := f(\tau x)$. We thus obtain a 
reflection positive Hilbert space $(\cH_K, \cH_K^+,\theta)$. 
We call such kernels $K$ {\it reflection positive} with respect to $(X,X_+, \tau)$. 

In this context, the space $\hat\cE$ can be identified with the reproducing kernel 
space $\cH^{K^\tau} \subeq \C^{X_+}$, where $q$ corresponds to the map 
\[ q \: \cE_+ \to \cH_{K^\tau}, \quad 
q(f)(x) := f(\tau(x)).\] 
In fact, the space $\hat\cE$ is generated by the elements 
$(q(K_x))_{x \in X_+}$, and we have 
\[ \la q(K_x), q(K_y) \ra = \la \theta K_x, K_y \ra = \la K_{\tau x}, K_y \ra 
= K(\tau x, y) = K^\tau(x,y).\] 
Accordingly, the function on $X_+$ corresponding to 
$q(f) \in \hat\cE$ is given by 
\[ x \mapsto \la q(K_x), q(f) \ra = \la \theta K_x, f \ra = \la K_{\tau x}, f \ra
= f(\tau x).\] 

(b) If $(\cE, \cE_+,\theta)$ is a reflection positive Hilbert space, 
then the scalar product defines a reflection positive kernel 
$K(v,w) := \la v, w \ra$ with respect to $(\cE,\cE_+,\theta)$. In this sense all 
reflection positive Hilbert spaces can be obtained in the context of (a), which 
provides a ``non-linear'' setting for reflection positive Hilbert spaces. 
\end{example}

\begin{defn} \mlabel{def:2.5} 
Let $\tau \: G \to G$ be an involutive automorphism of the group~$G$ 
and $G_+ \subeq G$ be a subset. 
A function $\vphi \colon G \to \C$ is 
called {\it reflection positive} with respect to $(G,G_+,\tau)$  if the kernel 
$K(x,y) := \vphi(xy^{-1})$ is reflection positive with respect to $(G,G_+,\tau)$ 
in the sense of Example~\ref{ex:1.3}(a). 
These are two simultaneous positivity conditions, namely that the kernel 
$\vphi(gh^{-1})_{g,h \in G}$ is positive definite on $G$ and that the kernel 
$\vphi(s\tau(t)^{-1})_{s,t \in G_+}$ is positive definite on $G_+$.  

Of particular importance is the case where $G_+ = S$ is a subsemigroup invariant 
under $s \mapsto s^\sharp := \tau(s)^{-1}$, so that $(S,\sharp)$ is an involutive 
semigroup and the positive definiteness of the kernel 
$K^\tau$ means that $\vphi$ is a positive definite function on 
$(S,\sharp)$. 
\end{defn}

In the following we 
write $\R_+ := [0,\infty)$ for the set of non-negative real numbers. 

\begin{ex}
(a) Prototypical examples are the functions 
$\vphi_\lambda(t) := e^{-\lambda|t|}$, $\lambda \geq 0$, for $(\R,\R_+, -\id_\R)$. 
For this triple   every continuous reflection positive 
function has an integral representation 
$\vphi(t)= \int_0^\infty e^{-\lambda|t|}\, d\mu(\lambda)$ for a positive Borel measure 
$\mu$ on $\R_+$ (see \cite[Cor.~3.3]{NO14}). 
In \cite{NO14} we also discuss generalizations of this concept to distributions 
and obtain integral representations for the case where $G$ is a more general 
abelian Lie group.

(b) There are also important examples not related to subsemigroups. 
For  $\beta > 0$, the corresponding circle group 
$G := \R/\beta \Z$ and the domain 
$G_+ := \big[0,\frac{\beta}{2}\big] + \beta \Z\subeq G$, the functions 
$\vphi \: G \to \C$ correspond to $\beta$-periodic functions on $\R$. 
Such a function is called {\it reflection 
positive} if the kernel 
$K(x,y) := \vphi(x-y)$ is reflection positive for $(G,G_+, \tau)$ and 
$\tau(g) = g^{-1}$ in the sense of Example~\ref{ex:1.3}(a). 
Typical examples are the $\beta$-periodic functions whose 
restriction to $[0,\beta]$ is given by 
$f_\lambda(t) := e^{-t\lambda} + e^{-(\beta - t)\lambda}$ for $\lambda \geq 0$ 
(see also Remark~\ref{rem:4.9} below and \cite{NO15}). 
\end{ex}

\begin{definition} \mlabel{def:repo-aff}
We call a continuous function $\psi \colon G \to \R$ {\it reflection negative} 
with respect to $(G,G_+,\tau)$ if 
$\psi$ is a negative definite function on $G$ 
and the kernel $\psi(st^\sharp)_{s,t \in G_+}$ is negative definite. 
\end{definition}

\begin{remark} \mlabel{rem:schoen} According to Schoenberg's Theorem for kernels 
\cite[Thm.~3.2.2]{BCR84}, 
a function $\psi \colon G \to \C$ is 
reflection negative if and only if, for every $h > 0$, the function 
$e^{-h\psi}$ is reflection positive in the sense of 
Definition~\ref{def:2.5}. 
\end{remark}

\section{Negative definite kernels on intervals} 
\mlabel{sec:7}

In this section we describe an integral representation 
of negative definite functions on general open intervals $(a,b) \subeq \R$. 
This extends the L\'evy--Khintchine formula for Bernstein 
functions, i.e., non-negative negative definite functions on~$(0,\infty)$. 

\begin{defn} Let $-\infty \leq a < b \leq \infty$. 
A function $f \: (a,b) \to \R$ is said to be 
\begin{itemize}
\item {\it positive definite} if the kernel $f\big(\frac{x+y}{2}\big)$ 
is positive definite, and 
\item {\it negative definite} if the kernel $f\big(\frac{x+y}{2}\big)$ 
is negative definite. 
\end{itemize}
\end{defn}

We first recall some classical results by Bernstein, Hamburger 
and Widder concerning an intrinsic characterization of Laplace transforms. 

\begin{thm} {\rm(Widder; \cite{Wi34}, \cite[Thm.~VI.21]{Wi46})}  \mlabel{thm:widder}
A function $\vphi\: (a,b) \to \R$ is positive definite if and only 
if there exists a positive measure $\mu$ on $\R$ such that 
\[ \vphi(t) = \cL(\mu)(t) := \int_\R e^{-\lambda t}\, d\mu(\lambda) 
\quad \mbox{ for }  \quad t \in (a,b).\] 
This implies in particular that $\vphi$ is analytic. 
\end{thm}

\begin{lem} \mlabel{lem:7.1} Every negative definite function 
$\psi \: (a,b) \to \R$ is analytic. 
\end{lem}

\begin{prf} Since $\vphi := e^{-\psi}$ is positive definite, 
it is analytic by Widder's Theorem~\ref{thm:widder}, so that 
\break $\psi = - \log \vphi$ is analytic as well.   
\end{prf}

\begin{thm} \mlabel{thm:8.5} 
Let $f \: (a,b) \to \R$ be an analytic function. 
\begin{itemize}
\item[\rm(a)] {\rm(Hamburger)} The following are equivalent: 
  \begin{itemize}
  \item[\rm(i)] $f$ is positive definite, i.e., 
$f = \cL(\mu)$ for a positive Borel measure $\mu$ on $\R$. 
  \item[\rm(ii)] For every $c \in (a,b)$,  
the kernel $(f^{(i+j)}(c))_{i,j \in \N_0}$ is positive definite. 
  \item[\rm(iii)] There exists a $c \in (a,b)$ 
for which the kernel $(f^{(i+j)}(c))_{i,j \in \N_0}$ is positive definite. 
  \end{itemize}
\item[\rm(b)] {\rm(Widder)} The following are equivalent: 
  \begin{itemize} 
  \item[\rm(i)] $f = \cL(\mu)$ for a positive Borel measure $\mu$ on $[0,\infty)$. 
  \item[\rm(ii)] For every $c \in (a,b)$,  
the kernels $(f^{(i+j)}(c))_{i,j \in \N_0}$ and 
$(- f^{(1+i+j)}(c))_{i,j \in \N_0}$ are positive definite. 
  \item[\rm(iii)] There exists a $c \in (a,b)$ 
for which the kernels $(f^{(i+j)}(c))_{i,j \in \N_0}$ and 
$(- f^{(1+i+j)}(c))_{i,j \in \N_0}$ are positive definite. 
  \item[\rm(iv)] $f$ and $-f'$ are positive definite. 
  \end{itemize}
\end{itemize}
\end{thm}

\begin{prf} (a) The implications (i) $\Rarrow$ (ii) $\Rarrow$ (iii) are trivial. 
For the remaining implication we refer to 
\cite[Lemma~3]{Wi34}, \cite[Thm.~VI.19c]{Wi46} or \cite{Ham20}. 

(b) Again, the implications (i) $\Rarrow$ (ii) $\Rarrow$ (iii) are trivial 
and (iii) $\Rarrow$ (i) follows from \cite[Thm.~VI.19b]{Wi46}. 
That (ii), (iii) are equivalent to (iv) now follows from~(a).
\end{prf}

\begin{defn} Let $-\infty \leq a < b \leq \infty$. 
  \begin{itemize}
  \item A smooth function $\psi \: (a,b) \to \R$ is called 
{\it completely monotone} if 
$(-1)^k \psi^{(k)} \geq 0$ for every $k \in \N = \{1,2,\ldots\}$. 
\item A smooth function $\psi \: (0,\infty) \to \R$ is called a {\it Bernstein function} 
if $\psi \geq 0$ and $\psi'$ is completely monotone. 
  \end{itemize}
\end{defn}

The following theorem provides a characterization of completely monotone 
functions on $(0,\infty)$:

\begin{thm} {\rm(Hausdorff--Bernstein--Widder)} \mlabel{thm:bern1}
For a function $\vphi \: (0,\infty) \to [0,\infty)$, the following are equivalent: 
\begin{itemize}
\item[\rm(i)] $\vphi$  is completely monotone.  
\item[\rm(ii)] $\vphi$ is a Laplace transform of a positive measure on 
$[0,\infty)$. 
\item[\rm(iii)] $\vphi$ is decreasing and 
positive definite on the $*$-semigroup $((0,\infty),\id)$. 
\end{itemize}
\end{thm}

\begin{prf} The equivalence of (i) and (ii) follows from 
\cite[Thms.~1.4]{SSV10}. %\cite[Thm.~IV.12b]{Wi46})} 
We further know from \cite[Thm.~III.1.19]{Ne00}, a positive definite function 
$\vphi$ on $(0,\infty)$ is decreasing if and only if the 
corresponding representation is a representation by contractions. 
Therefore the equivalence of (ii) and (iii) is Nussbaum's 
Theorem \cite[Cor.~VI.2.11]{Ne00}. 
\end{prf}

\begin{rem} That an analytic function $\vphi \: (a,b) \to \R$ 
is a Laplace transform implies in particular
that $\vphi^{(2j)} \geq 0$ for each $j \in \N_0$, but this takes 
only the diagonal entries of the matrices in Theorem~\ref{thm:8.5}(a)(ii) 
into account. Likewise, condition (b)(ii) in 
Theorem~\ref{thm:8.5} implies that $f$ is completely monotone 
and Theorem~\ref{thm:bern1} shows that, on $(a,b) = (0,\infty)$, 
the stronger conditions in Theorem~\ref{thm:8.5} follow from complete monotonicity.
\end{rem}

The following lemma prepares our characterization of 
negative definite functions in Theorem~\ref{thm:L-K-interval} below. 

\begin{lem} \mlabel{lem:7.5} If $\psi$ is negative definite, then 
$-\psi''$ is positive definite, resp., a Laplace transform of a 
positive measure on~$\R$.   
\end{lem}

\begin{prf} In view of the analyticity of the function $- \psi''$ 
(Lemma~\ref{lem:7.1}), 
it suffices to show that, in some $c \in (a,b)$, the kernel 
$(-\psi^{(i+j)}(c))_{i,j \geq 1}$ is positive definite 
(Hamburger's Theorem~\ref{thm:8.5}(a)). This is equivalent 
to the matrices $(-\psi^{(i+j)}(c))_{1 \leq i,j \leq n}$ being positive semi-definite 
for every $n \in \N$. 

For $\delta \in \R$, let $(\Delta_\delta f)(t) := 
f(t) - f(t + \delta)$. 
Then 
\[ \lim_{\delta \to 0} \delta^{-j} (\Delta_\delta^j f)(t) = (-1)^j 
f^{(j)}(t) \] 
for any smooth function $f \: (a,b) \to \R$. 

Since $\psi$ is negative definite, for every $h > 0$, the function 
$\vphi_h := e^{-h\psi}$ is positive definite (Schoenberg's Theorem). 
Therefore 
\begin{equation}
  \label{eq:approx}
  \psi 
= \lim_{h \to 0+} \frac{1-e^{-h\psi}}{h} 
\end{equation}
shows that $\psi$ is a pointwise limit of positive multiples 
of functions of the form $1 - \vphi$, where $\vphi$ is positive 
definite. 

If $\vphi$ is positive definite, 
$c \in (a,b)$ and $c + 2n \delta \in (a,b)$, then the matrix 
$\big( (\Delta^{i+j}_\delta \vphi)(c)\big)_{0 \leq i,j \leq n}$ 
is positive semidefinite 
(cf.\ the proof of the theorem in \cite{Wi34}). 
This implies that, for any negative definite function of the form 
$\tilde\psi = 1 - \vphi,$ the matrix 
\[ ( -\big(\Delta^{i+j}_\delta \tilde\psi)(c)\big)_{1 \leq i,j \leq n} 
=  \big( (\Delta^{i+j}_\delta \vphi)(c)\big)_{1 \leq i,j \leq n} \] 
is positive semidefinite. Since any negative definite function 
$\psi$ is a pointwise 
limit of positive multiples of such functions $\tilde\psi$, we see that 
$\big( -(\Delta^{i+j}_\delta \psi)(c)\big)_{1 \leq i,j \leq n}$
is positive semidefinite. For $\delta \to 0$, we thus obtain that 
the matrix 
$\big(-(-1)^{i+j} \psi^{(i+j)}(c)\big)_{1 \leq i,j \leq n}$ 
is positive semidefinite, and this implies that 
$\big( -\psi^{(i+j)}(c)\big)_{1 \leq i,j \leq n}$
is positive semidefinite as well. This completes the proof.
\end{prf}

The following theorem is a slight generalization of \cite[Thm.~4.2]{Ho67} 
which only deals with intervals of the form $[0,b]$. Our overall strategy 
is similarly to Horn's. 

\begin{thm}{\rm(L\'evy--Khintchine formula for open intervals)} 
  \mlabel{thm:L-K-interval}
Fix $t_0 \in (a,b)$.  Then $\psi \: (a,b) \to \R$ is negative definite 
if and only if there exists a 
positive measure $\mu$ on $\R$ and $c,d \in \R$ such that 
\begin{equation}
  \label{eq:leki}
\psi(t) = c + d(t-t_0) + \int_{-\infty}^\infty e_\lambda(t)e^{-\lambda t_0}
\, d\mu(\lambda),
\end{equation}
where 
\begin{equation}
  \label{eq:e-lambda}
e_\lambda(t) := \begin{cases}
\frac{1 - \lambda (t-t_0) - e^{-\lambda (t-t_0)}}{\lambda^2} 
& \text{ for } \lambda\not=0 \\ 
-(t-t_0)^2/2  & \text{ for } \lambda=0.  
\end{cases}
\end{equation}
Then $c = \psi(t_0)$, $d = \psi'(t_0)$ 
and $\mu$ is uniquely determined by $-\psi'' = \cL(\mu)$. 
A positive measure $\mu$ 
occurs for some~$\psi$ 
if and only if its Laplace transform $\cL(\mu)$ is finite on $(a,b)$. 
\end{thm}

\begin{prf} ``$\Rarrow$'': Let $\psi$ be negative definite. 
We prove the existence of the integral representation. 
{}From Lemma~\ref{lem:7.5} we know that $-\psi''$ is positive definite. 
Therefore Widder's Theorem~\ref{thm:widder} implies the existence of a positive 
measure $\mu$ on $\R$ with $\psi''(t) = - \cL(\mu)(t)$ for every $t \in (a,b)$. 
Since $\psi$ is a smooth function, we have for $t \in (a,b)$ the formula 
\begin{equation}
  \label{eq:taylor}
\psi(t) = \psi(t_0) + (t-t_0) \psi'(t_0) 
+ \int_{t_0}^t \psi''(s) \cdot(t-s)\, ds.
\end{equation}
We thus put $c := \psi(t_0)$ and $d := \psi'(t_0)$. For the 
third term in \eqref{eq:taylor} we find with Fubini's Theorem: 
\begin{equation}
  \label{eq:19}
\int_{t_0}^t \psi''(s) (t-s)\, ds
= -\int_{t_0}^t \int_\R e^{-\lambda s}\, d\mu(\lambda) \cdot (t-s)\, ds
= \int_\R \Big(\int_{t_0}^t e^{-\lambda s}(s-t)\, ds\Big)\, d\mu(\lambda).
\end{equation}
The statement now follows by integration by parts. But we would also
like to point out another argument.
As the functions $e_\lambda$ in \eqref{eq:e-lambda} are uniquely determined 
by 
\[ e_\lambda(t_0) = e_\lambda'(t_0) = 0 \quad \mbox{ and }\quad 
e_\lambda''(t) = - e^{-\lambda(t-t_0)} = - e^{\lambda t_0} e^{-\lambda t},\] 
it follows from \eqref{eq:taylor} that 
\begin{equation}
  \label{eq:20}
e^{-\lambda t_0} e_\lambda(t) = \int_{t_0}^t e^{-\lambda s}(s-t)\, ds.
\end{equation}
Combining \eqref{eq:taylor}, \eqref{eq:19} and \eqref{eq:20}, we thus obtain 
the stated integral representation from 
\[ \psi(t) = \psi(t_0) + (t-t_0) \psi'(t_0) 
+ \int_\R e_\lambda(t) e^{-\lambda t_0}\, d\mu(\lambda).\]

``$\Larrow$'': To see that all functions with such an integral representation 
are negative definite, it suffices to observe that all affine functions 
are negative definite and that all functions $e_\lambda$ are negative 
definite. For $\lambda \not=0$, this follows from the negative 
definiteness of the functions 
$1 - e^{-\lambda t}$ (because $e^{-\lambda t}$ is positive definite) 
and of all affine functions. For $\lambda =0$ it follows from 
$e_0 = \lim_{\lambda \to 0} e_\lambda$. \\

Next we show that $\psi$ determines the constants 
$c,d$ and the measure $\mu$. If $\psi$ satisfies \eqref{eq:leki}, then, by reversing our previous
arguments to show that (\ref{eq:19}) implies (\ref{eq:leki}), and using that $e_\lambda (t)e^{-\lambda t}$ 
is positive, Fubini's Theorem implies that
\[ \int_{-\infty}^\infty e_\lambda(t)e^{-\lambda t_0} \, d\mu(\lambda)  
=  \int_{t_0}^t  \int_{-\infty}^\infty e^{-\lambda s}\, d\mu(\lambda)\cdot (s-t) \, ds.\]
In particular, the existence of
the integrals in \eqref{eq:leki} implies that 
$\cL(\mu)(s)$ is finite for every $s$ in a subset of full measure in $(a,b)$. 
Therefore the convexity of $\cL(\mu)$ shows that 
$\cL(\mu)(t) < \infty$ for every $t \in (a,b)$. This leads to 
\[ \psi(t) = c + d(t-t_0) - \int_{t_0}^t \cL(\mu)(s)\cdot (t-s)\, ds, \] 
and this entails that $c = \psi(t_0), d = \psi'(t_0)$ and 
$\psi''(t) = - \cL(\mu)(t)$ for every $t \in (a,b)$. 
Therefore $c,d$ and $\mu$ are uniquely determined by~$\psi$. 
\end{prf}

\begin{rem} \mlabel{rem:8.8} 
(a) Theorem~\ref{thm:widder} has been generalized by Shucker 
in \cite[Thm.~5]{Sh84} to continuous positive definite 
functions $r \: \Omega \to \C$
on convex domains in a real vector space 
(see also \cite[Thm.~4.11]{NO14} for a related integral representation 
of $1$-bounded distributions on cones). 
For intervals in $\R^n$, there are corresponding 
results by Widder and Akhiezer \cite[Thm.~6.5.12]{BCR84}. 

(b) For integral representations of negative definite functions 
on $\Q_+$, see \cite[Prop.~6.5.13]{BCR84}. 

(c) According to \cite[Thm.~6.5.14]{BCR84}, continuous negative definite 
functions on the closed half line $[0,\infty)$ are of the form 
\begin{equation}
  \label{eq:LK-BCR}
\psi(t) = c + dt- f t^2 + 
\int_{\R^\times} \Big(1 - e^{\lambda t} + \frac{\lambda t}{1 + \lambda^2}\Big)\, 
d\mu(\lambda),
\end{equation}
where $a,b \in  \R$, $c \geq 0$ and $\mu$ is a positive measure on $\R^\times$ 
satisfying 
\begin{equation}
  \label{eq:21}
\int_{0 < |\lambda| \leq 1} \lambda^2\, d\mu(\lambda) < \infty \quad \mbox{ and } \quad 
\int_{|\lambda| >1} e^{\lambda t}\, d\mu(\lambda) < \infty 
\quad \mbox{ for }\quad t\geq 0.
\end{equation}
This matches the integral representation in Theorem~\ref{thm:L-K-interval} 
for $a = 0$ and $b = \infty$ because \eqref{eq:LK-BCR}  
implies that 
\[ -\psi''(t) =  2 f  + \int_{\R^\times} e^{\lambda t} \lambda^2\, d\mu(\lambda),\] 
so that the measure in  Theorem~\ref{thm:L-K-interval} is 
$2f \delta_0 + \lambda^2 \, d\mu(-\lambda).$ 
Here conditions \eqref{eq:21} on the measure $\mu$ correspond to the additional 
requirement that $\psi(0+) = \lim_{t \to 0+} \psi(t) < \infty$. 

(d) In \cite[Ex.~8.1.11]{BCR84} it is shown that 
a negative definite function $\psi \: (0,\infty) \to \R$ which is bounded 
from below on all intervals $[c,\infty)$ has a representation 
\begin{equation}
  \label{eq:22}
\psi(s) - \psi(t) = \int_{\R^\times_+} (e^{-\lambda t} - e^{-\lambda s})\, d\mu(\lambda) 
+ \alpha \cdot (s-t)\quad \mbox{ for }\quad s,t > 0,
\end{equation}
where 
$\alpha \geq 0$ and $\mu$ is a positive measure on $(0,\infty)$. 
For $(s,t) = (t,1)$, this leads to 
\begin{equation}
  \label{eq:23}
\psi(t) = \psi(1) +  \alpha \cdot (t-1) + 
\int_{\R^\times_+} (e^{-\lambda} - e^{-\lambda t})\, d\mu(\lambda) 
\quad \mbox{ for }\quad t > 0,
\end{equation}
which is a special case of the integral representation 
in Theorem~\ref{thm:L-K-interval} which immediately implies \eqref{eq:22} 
by subtraction. 

A typical example is the function 
\begin{equation}
  \label{eq:log-int}
\log t = \int_0^\infty (e^{-\lambda} - e^{-\lambda t})\frac{d\lambda}{\lambda} 
\quad \mbox{ for } \quad t > 0
\end{equation}
(see also Examples~\ref{ex:bernstein}(c)). 
\end{rem}

\section{Increasing negative definite functions on $(0,\infty )$} 

Now we turn to the special case where $a = 0$ and $b = \infty$. 
Then $S := ((0,\infty),+)$  is a $*$-semigroup with respect to 
$s^* = s$ for $s \in S$. 
Theorem~\ref{thm:L-K-interval} provides in particular 
an integral formula for an arbitrary 
negative definite function on $(0,\infty)$. For the applications 
in representation theory, we are also interested in the increasing 
negative definite functions which are not necessarily non-negative, 
resp., Bernstein functions. These are characterized in
 Bernstein's classical theorem: 

\begin{thm} {\rm(Bernstein)} {\rm(\cite[Thms.~1.4, 3.2, 3.6]{SSV10})} 
\mlabel{thm:lk-bernstein2}
For a function $\psi \: (0,\infty) \to [0,\infty)$, the following are equivalent: 
\begin{itemize}
\item[\rm(i)] $\psi$  is a Bernstein function. 
\item[\rm(ii)] For every $h> 0$, the function 
$e^{-h \psi}$ is completely monotone. In particular, $\psi$ is negative definite 
and increasing. 
\item[\rm(iii)] There exist $a, b \geq 0$ and a positive measure 
$\sigma$ on $(0,\infty)$ with $\int_0^\infty (1 \wedge \lambda)\, d\sigma(\lambda) < \infty$ 
such that 
\[ \psi(t) = a + bt + \int_0^\infty (1 - e^{-\lambda t})\, d\sigma(\lambda) \] 
(L\'evy--Khintchine representation). Then 
$a = \lim_{t \to 0} \psi(t),  b = \lim_{t \to \infty} \frac{\psi(t)}{t},$
and $\sigma$ is also uniquely determined by~$\psi$.
\end{itemize}
\end{thm}

\begin{rem} (a) Theorem~\ref{thm:lk-bernstein2} 
describes the non-negative real-valued 
negative definite functions on the $*$-semigroup $S := ((0,\infty),\id)$, 
but not every negative definite function on $S$ is real-valued and non-negative.
For instance, all affine functions $\psi(t) = a + bt$, $a,b \in \R$, 
are negative definite (Theorem~\ref{thm:L-K-interval}). 

(b) For any $c > 0$, the function $\psi(t) = 1 - c e^{-\lambda t}$ 
is negative definite on $(\R_+, \id)$ because the corresponding kernel 
is 
$\psi(s+t) = 1 - c e^{-\lambda s} e^{-\lambda t},$ 
where both summands are negative definite. 
However, only for $c \leq 1$, we obtain a Bernstein function on $\R_+$. 
In this case 
\[ \psi(t) = (1 - c) + c(1 - e^{-\lambda t}) \] 
is the corresponding L\'evy--Khintchine representation. 

(c) If $\vphi \: (0,\infty) \to \R$ is a non-zero decreasing positive definite function on 
$S$, then $\vphi(t) > 0$ for every $t > 0$, so that 
$\psi := - \log \vphi \: (0,\infty) \to \R$ is increasing, 
but in general $\psi$ may take negative values. 

If $\psi$ is a Bernstein function, then $\vphi = e^{-\psi}\leq 1$. 
If, conversely, 
$\vphi \leq C$ is bounded, then $\psi \geq - \log C$, so that $\psi + \log C \geq 0$. 
If, more generally, there exist 
$a,b \in \R$ with $\vphi(t) \leq e^{a + bt}$ for all $t > 0$, then 
$\psi - a - bt \geq 0$. 
\end{rem}

\begin{exs} (of Bernstein functions)  \mlabel{ex:bernstein} 

(a) For $\alpha \in \R$, the function $\psi(t) =  t^\alpha$ on $(0,\infty)$ 
has positive values. From $\psi'(t) = \alpha t^{\alpha-1}$ we derive that 
$\alpha \geq 0$ is necessary for $\psi$ to be a Bernstein function. 
Further, $\psi''(t) = \alpha(\alpha-1) t^{\alpha-2}$ shows that 
$\alpha \leq 1$ is also necessary. 

Conversely, $\alpha = 0,1$ lead to the Bernstein functions $1$ and $t$, 
and for $0 < \alpha < 1$, the function $t^\alpha$ is also Bernstein 
because 
$(-1)^{k-1} \psi^{(k)}(t) = \alpha(1-\alpha) \cdots (k - 1 - \alpha) 
t^{\alpha -k} \geq 0$ 
for every $t >0$. Its L\'evy--Khintchine representation is 
\[ t^\alpha = \frac{\alpha}{\Gamma(1-\alpha)} \int_0^\infty (1- e^{-\lambda t}) 
\lambda^{-1-\alpha}\, d\lambda \quad \mbox{ for } \quad 0 < \alpha < 1, t > 0\] 
(see the proof of \cite[Cor.~3.2.10]{BCR84}). 
 
(b)  $\frac{t}{1+t} = \int_0^\infty (1 - e^{-\lambda t}) e^{-\lambda}\, d\lambda$ for $t > 0$. 

(c) For $\alpha > 0$, we have 
\[ t^{-\alpha} = \frac{1}{\Gamma(\alpha)} \int_0^\infty \lambda^{\alpha-1} e^{-\lambda t}\, d\lambda 
= \cL(\mu_\alpha)(t) \quad \mbox{ for }  \quad 
d\mu_\alpha(\lambda) = \lambda^{\alpha-1}\, d\lambda.\] 
This implies that all positive powers of the function 
$\vphi(t) = t^{-1}$ are positive definite on the semigroup~$\R_+$ 
(Theorem~\ref{thm:bern1}), so that 
\[ \psi(t) := \log t  \] 
is a negative definite function on the additive semigroup $(0,\infty)$ by 
Schoenberg's Theorem (Remark~\ref{rem:schoen}). 
We now derive an integral representation for the function $\psi$ 
which is not bounded from below and in particular not Bernstein. 
Nevertheless, the shifted function $\log(1 + t)$ is Bernstein 
(see \eqref{eq:log}). 
Starting with 
\[ 1 - e^{-\alpha \log t} = 1 - t^{-\alpha} = 
\frac{1}{\Gamma(\alpha)} \int_0^\infty \lambda^{\alpha-1} e^{-\lambda} (1-e^{-\lambda(t-1)})
\, d\lambda,\]
we consider with Lebesgue's Dominated Convergence Theorem 
\begin{align*}
\log t
&= \lim_{\alpha \to 0} \frac{1 - e^{-\alpha \log t}}{\alpha}
= \lim_{\alpha \to 0} \frac{1}{\alpha\Gamma(\alpha)} 
\int_0^\infty \lambda^{\alpha-1} e^{-\lambda} (1-e^{-\lambda(t-1)})\, d\lambda\\
&= \lim_{\alpha \to 0} \frac{1}{\Gamma(\alpha+1)} 
\int_0^\infty \lambda^{\alpha-1} e^{-\lambda} (1-e^{-\lambda(t-1)})\, d\lambda 
= \int_0^\infty  (1-e^{-\lambda(t-1)})\, e^{-\lambda}\frac{d\lambda}{\lambda} \\
&= \int_0^\infty (e^{-\lambda} -e^{-\lambda t})\, \frac{d\lambda}{\lambda}. 
\end{align*}
This leads to the integral representation 
\begin{equation}
  \label{eq:log}
\log(1 + t) = \int_0^\infty (1 - e^{-\lambda t}) e^{-\lambda}\, \frac{d\lambda}{\lambda} 
\quad \mbox{ for } \quad t > -1 
\end{equation}
which exhibits $\log(1 + t)$ as a Bernstein function on $(0,\infty)$.
\end{exs}

\begin{ex} \mlabel{ex:power}
(a)  (\cite[Ex.~6.5.15]{BCR84}) 
For $0 < \alpha \leq 2$, the function 
\[ \psi_\alpha \: [0,\infty) \to \R, \quad \psi_\alpha(t) :=
\begin{cases}
  t^\alpha & \text{ for } 0 < \alpha \leq 1 \\ 
  -t^\alpha & \text{ for } 1 \leq \alpha \leq 2
\end{cases}
\]
is negative definite on the semigroup $(\R_+,\id)$. 

(b) (\cite[Ex.~6.5.16]{BCR84})  
The function $\psi(t) = - t \log t$ is negative definite on $[0,\infty)$, 
i.e., for every $h > 0$, the function $e^{-h \psi(t)} = e^{h t \log t} = t^{ht}$ 
is positive definite. 
\end{ex}

We now turn to our generalization of Bernstein's Theorem characterizing 
the increasing negative definite functions on $(0,\infty)$. 

\begin{rem} \mlabel{rem:7.9} For a $C^k$-function, we consider the difference operators 
\[ (\Delta_\delta f)(t) :=f(t) - f(t + \delta) 
= (-1) \int_t^{t+\delta} f'(s)\, ds.\] 
An easy induction shows that 
\[ (\Delta_h^k f)(t) = (-1)^k 
\int_t^{t+h} \int_{t_1}^{t_1 + h} \cdots \int_{t_{k-1}}^{t_{k-1}+h} 
f^{(k)}(t_k)\, dt_k \cdots dt_1.\] 
Therefore $(\Delta_h^k f)(t) \geq 0$ for all $t$ and sufficiently 
small $h$ (depending on $t$) is equivalent to the condition $(-1)^k f^{(k)} \geq 0$. 
\end{rem}

\begin{lem} \mlabel{lem:7.10} If $\psi \: (0,\infty) \to \R$ 
is negative definite and increasing, then $\psi'$ is completely monotone.
\end{lem}

\begin{prf} Lemma~\ref{lem:7.1} implies that $\psi$ is smooth, so that 
we have to show that $(-1)^k \psi^{(k+1)} \geq 0$ for $k \in \N_0$. 
In view of Remark~\ref{rem:7.9}, it suffices to show that 
$\Delta_\delta^k \psi' \geq 0$ for $\delta> 0$. This in turn will follow if 
$-\Delta_\delta^{k+1} \psi \leq 0$ for $k \in \N_0$. 

If $\psi = 1 -\vphi$, where $\vphi$ is positive definite, then 
$\vphi$ is decreasing, hence a Laplace transform of a measure on 
$[0,\infty)$ and therefore completely monotone  
(Theorem~\ref{thm:lk-bernstein2}). This implies that 
\[ \Delta_\delta^{k+1} \psi  =  -\Delta_\delta^{k+1} \vphi  \leq 0.\] 
Next we use $\psi = \lim_{h \to 0} \psi_h$ with $\psi_h = \frac{1 - e^{-h\psi}}{h}$ 
to see that 
$\Delta_\delta^{k+1} \psi  =  \lim_{h \to 0} \Delta_\delta^{k+1} \psi_h \leq 0.$ 
This completes the proof. 
\end{prf}

\begin{thm}{\rm(L\'evy--Khintchine formula for increasing 
negative definite functions on $\R_+$)} 
  \mlabel{thm:L-K-incr}
A function $\psi \: (0, \infty) \to \R$ is negative definite and increasing 
if and only if there exists a positive measure $\mu$ on $\R$ and $c\in \R$ such that 
\begin{equation}
  \label{eq:leki-inc}
\psi(t) = c + \int_{[0,\infty)} f_\lambda(t) \, d\mu(\lambda),
\quad \mbox{ where } \quad 
f_\lambda(t) := 
\begin{cases}
\frac{e^{-\lambda} - e^{-\lambda t}}{\lambda} 
& \text{ for } \lambda\not=0 \\ 
t-1  & \text{ for } \lambda=0.  
\end{cases}
\end{equation}
Then $c = \psi(1)$ and $\mu$ is uniquely determined by $\psi' = \cL(\mu)$. 
A positive measure $\mu$ 
occurs if and only if its Laplace transform $\cL(\mu)$ is finite on $(0,\infty)$. 
\end{thm}

\begin{prf} This result could also be derived from 
\eqref{eq:23} in Remark~\ref{rem:8.8}(d), but we also give an independent proof. 
We use a similar argument as in the proof of 
Theorem~\ref{thm:L-K-interval}, but here it is simpler. 
For the  existence of the integral representation, 
we first use Lemma~\ref{lem:7.10} to see that 
$\psi'$ is completely monotone, hence of the form 
$\cL(\mu)$ for a measure $\mu$ on $[0,\infty)$. 
We put $c := \psi(1)$ and observe that 
\begin{equation}
  \label{eq:taylor2}
\psi(t) = \psi(1) + \int_{1}^t \psi'(s)\, ds 
= \psi(1) + \int_{1}^t \int_0^\infty e^{-\lambda s}\, d\mu(\lambda)\, ds 
= \psi(1) + \int_0^\infty \Big(\int_{1}^t  e^{-\lambda s}\, ds\Big)\, d\mu(\lambda). 
\end{equation}
The functions $f_\lambda$ are characterized by $f_\lambda(1) = 0$ and 
$f_\lambda'(t) = e^{-\lambda t}$. We thus obtain the desired 
integral representation of~$\psi$.

To see that all functions with such an integral representation 
are negative definite and increasing, it suffices to observe that 
the functions $f_\lambda$ all have this property. 

To see that $\psi$ determines the measure $\mu$ uniquely, 
one argues as in the proof of Theorem~\ref{thm:L-K-interval} that 
the integral representation for $\psi$ implies that 
$\cL(\mu)$ is finite on $(0,\infty)$ and coincide with $\psi'$. 
\end{prf}

\begin{rem} (a) From \eqref{eq:leki-inc} we obtain 
\[ \psi(t) = c  + \mu(\{0\}) (t-1) + \int_{(0,\infty)} 
(e^{-\lambda} - e^{-\lambda t})\, \frac{d\mu(\lambda)}{\lambda},\] 
which is \eqref{eq:23} in Remark~\ref{rem:8.8}(d). 

(b) The increasing function $\psi$ in \eqref{eq:leki-inc} 
is non-negative, i.e., a Bernstein function, if and only if 
\[ \psi(0+) = \lim_{t \to 0} \psi(t) 
= c + \int_{[0,\infty)} f_\lambda(0)\, d\mu(\lambda)
= c - \mu(\{0\}) + \int_0^\infty \frac{e^{-\lambda} - 1}{\lambda}\, d\mu(\lambda) \] 
exists and is non-negative. With (a) we then obtain 
\begin{align*}
 \psi(t) 
&= c  + \mu(\{0\}) (t-1) + \int_{(0,\infty)}  
(e^{-\lambda} - e^{-\lambda t})\, \frac{d\mu(\lambda)}{\lambda} \\ 
&= \psi(0+)  + \mu(\{0\}) t + \int_{(0,\infty)}  
(1 - e^{-\lambda t})\, \frac{d\mu(\lambda)}{\lambda}.  
\end{align*}
This corresponds to the L\'evy--Khintchine 
formula for Bernstein functions in Theorem~\ref{thm:lk-bernstein2}. 
\end{rem}

\section{Reflection positive functions on intervals} 
\mlabel{sec:5}

In this section we eventually turn to reflection positive functions 
on intervals. We consider the interval $(-a,a)$ for some $a > 0$ 
and the reflection $\tau(t) = -t$ about the midpoint. 

\begin{defn} We call a function 
$\vphi \: (-a,a) \to \R$ {\it reflection positive} 
if both kernels 
\[ \vphi\Big(\frac{t-s}{2}\Big)_{-a < s,t < a} \quad \mbox{ and } \quad 
\vphi\Big(\frac{t+s}{2}\Big)_{0< s,t < a} \] 
are positive definite. 
\end{defn}

This corresponds to the 
situation of Example~\ref{ex:1.3}(a) where $X = (-a,a)$, $X_+ = (0,a)$ 
and $\tau(x) = -x$. 

Reflection positivity implies that $\vphi(-t) = \vphi(t) 
= \oline{\vphi(t)}$. 
By Widder's Theorem~\ref{thm:widder}, 
there exists a positive measure $\mu$ on $\R$ with 
\begin{equation}
  \label{eq:lapl-abs}
\vphi(t) = \cL(\mu)(|t|) = \int_\R e^{-\lambda |t|}\, d\mu(\lambda) \quad \mbox{ 
for } \quad |t| < a.
\end{equation}
For all these functions the kernel 
$\vphi\big(\frac{t+s}{2}\big)$ is positive definite on $(0,a)$. Therefore 
$\vphi$ is reflection positive if and only if the kernel 
$\vphi\big(\frac{t-s}{2}\big)$ is positive definite on $(-a,a)$. 
Here the main point is to relate this condition to properties 
of the measure~$\mu$.

\begin{ex} (a) 
For $\lambda \geq 0$, the functions 
$\vphi_\lambda(t) := e^{-\lambda |t|}$ 
(multiples of euclidean Green's functions \cite{DG13}) are positive definite on 
$\R$. Therefore $\cL(\mu)(|t|)$ is reflection positive if $\mu$ is 
supported by $[0,\infty)$.   

(b) (cf.\ \cite[Ex.~2.3]{NO15}) Basic examples of  positive definite 
$\beta$-periodic functions on $\R$ are given by  the functions 
$f_\lambda$ satisfying 
\[ f_\lambda(t) = e^{-t\lambda} + e^{-(\beta - t)\lambda} 
= 2 e^{-\beta \lambda/2} \cosh((\textstyle{\frac{\beta}{2}}-t)\lambda)
\quad \mbox{ for } \quad 0 \leq t \leq \beta, \lambda \geq 0\] 
(multiples of thermal euclidean Green's functions \cite{DG13}). 
For $|t| < \beta$, we have
\begin{equation}
  \label{eq:f-lambda}
f_\lambda(t) = f_\lambda(|t|) 
= e^{-|t|\lambda} + e^{-(\beta - |t|)\lambda} 
= e^{-|t|\lambda} + e^{-\beta\lambda} e^{|t|\lambda}.
\end{equation}
Hence, for reflection positivity, it is not necessary that 
$\mu$ is supported by the positive half line, as in (a). 

Given $a > 0$ and a positive measure $\mu$ on 
$[0,\infty) \times [a,\infty)$, it follows that 
the function 
\begin{equation}
  \label{eq:doubleint}
f(t) := \int_0^\infty \int_a^\infty 
e^{-\lambda|t|} + e^{-\beta\lambda} e^{\lambda|t|}\, d\mu(\lambda,\beta) 
\end{equation}
is reflection positive on $(-a,a)$ whenever the integrals are finite. 
\end{ex} 

\begin{rem} For a function $\phi = \cL(\mu)$ as in \eqref{eq:lapl-abs}, 
a necessary condition for positive definiteness on $(-a,a)$ is that 
$\vphi(t) \leq \vphi(0)$ for $|t| < a$ because the positive definiteness 
of the matrix 
$\pmat{\vphi(0) & \vphi(t) \\ \vphi(-t) & \vphi(0)}$ implies 
$|\vphi(t)|^2 = \vphi(t)\vphi(-t) \leq \vphi(0)^2$. 
If, in addition, $\vphi$ is reflection positive, then 
$\vphi(t) = \cL(\mu)(|t|)$ and we obtain the condition  
$\vphi(t) \leq \mu(\R)$ for $0 \leq t < a$ and by convexity of the 
function $\vphi$ on $(0,a)$ we get
\begin{equation}
  \label{eq:ineq}
\vphi(a_-) = \lim_{t \to a-} \vphi(t) 
= \int_\R e^{- \lambda a}\, d\mu(\lambda) \leq \vphi(0),\quad \mbox{ resp.,} 
\quad 
\int_\R (1- e^{- \lambda a})\, d\mu(\lambda) \geq 0.
\end{equation}
For $\mu = \delta_{\lambda_0} + c \delta_{-\lambda_0}$,  condition \eqref{eq:ineq} means that $(1-e^{-\lambda_0 a})  + c (1 - e^{\lambda_0 a}) \geq 0,$ which is 
\[ 1 - e^{-\lambda_0 a} \geq c(e^{\lambda_0 a}-1) = e^{\lambda_0 a} c(1- e^{-\lambda_0 a}) 
\quad \mbox{ resp.} \quad 
c \leq e^{-\lambda_0 a}.\] 
Note that the maximal value of $c$ is precisely the constant showing 
up in \eqref{eq:f-lambda} with $a=\beta$.
\end{rem}

We now use P\'olya's Theorem to obtain sufficient 
conditions for positive definiteness on some interval $(-a,a)$: 

\begin{thm} {\rm(P\'olya)} \mlabel{thm:polya} {\rm(\cite[Thm.~4.3.1]{Luk70})}
If $\vphi \: \R \to [0,\infty)$ is an even continuous function 
which is convex on $[0,\infty)$ and satisfies 
$\lim_{t \to \infty} \vphi(t) = 0$, then 
$\vphi$ is positive definite.   
\end{thm}

\begin{cor} \mlabel{cor:polya2} 
If $\vphi \: \R \to [0,\infty)$ is an even continuous function 
which is convex and decreasing on $[0,\infty)$, 
then $\vphi$ is positive definite.   
\end{cor}

\begin{prf} Since $\vphi$ is decreasing, the limit 
$c := \lim_{t \to \infty} \vphi(t) \geq 0$ exists. Now P\'olya's Theorem 
implies that $\vphi - c$ is positive definite, and since the constant 
function $c$ is also positive definite, the assertion follows. 
\end{prf}

\begin{lem}\mlabel{le:5.6} Let $a > 0$ and $\psi \: [0,a] \to [0,\infty)$ be a convex function. 
Then $\psi$ extends to a non-negative decreasing convex function 
$\vphi \: [0,\infty) \to [0,\infty)$ if and only if $\psi'(a-) \leq 0$. 
\end{lem}

\begin{prf} If $\psi$ extends to a decreasing convex 
function $\vphi \: [0,\infty) \to [0,\infty)$, then 
$\psi'(a-) = \vphi'(a-) \leq 0$. 
If, conversely, $\psi'(a-) \leq 0$, then 
we extend $\psi$ by the constant function $t \mapsto \psi(a)$ on $(a,\infty)$. 
Then $\vphi$ is a decreasing convex function extending~$\psi$. 
\end{prf}

Combining the preceding lemma with Corollary~\ref{cor:polya2}, we obtain: 
\begin{prop} \mlabel{prop:polya} 
Let $a > 0$ and $\psi \: [0,a] \to [0,\infty)$ be a convex function. 
If $\psi'(a-) \leq 0$, then the function $\vphi(t) := \psi(|t|)$ on $[-a,a]$ 
is positive definite in the sense that the kernel 
$\vphi\big(\frac{t-s}{2}\big)_{|t|,|s| \leq a}$ is positive definite. 
\end{prop}

\begin{thm} \mlabel{thm:5.9} 
{\rm(Characterization of reflection positive functions on $[-a,a]$)} 
Let $\mu$ be a finite positive measure on $\R$ for which 
$\vphi(t) := \cL(\mu)(|t|)$ exists for $|t| \leq a$. 
\begin{itemize}
\item[\rm(a)] If $\cL(\mu)'(a-) \leq 0$, then $\vphi$ is reflection positive on $[-a,a]$ and extends to a symmetric positive definite function on $\R$. 
\item[\rm(b)] If $\vphi$ is reflection positive on $[-a,a]$ and non-constant, 
then there exists an element $b \in (0,a]$ with $\cL(\mu)'(b-) < 0$. 
\end{itemize}
\end{thm}

\begin{prf} If $\cL(\mu)'(a-) \leq 0$, then Proposition~\ref{prop:polya} 
implies that $\vphi$ is reflection positive on $[-a,a]$. The extension to $\R$ follows from
Lemma \ref{le:5.6}.

Suppose, conversely, that this is the case, 
and that $\vphi$ is not constant. Then there exists a $t \in (0,a]$ 
with $\vphi(t) < \vphi(0)$ because $\vphi(t) \leq \vphi(0)$ follows from positive 
definiteness of the kernel $\vphi\big(\frac{t-s}{2}\big)$. 
This means that $\vphi$ is not increasing on $[0,t]$, so that there exists a 
$t' \in (0,t)$ with $\vphi'(t'-) < 0$. 
\end{prf}

\begin{rem} If $b$ is as in (2) then it follows from part (a) that the positive definite function $\vphi|_{[-b,b]}$ extends to a positive definite function
on $\R$. But this extension does not have to agree with $\vphi$ on the interval $]b,a]$ if $b<a$.
\end{rem}

\begin{rem} (The $\beta$-periodic case) \mlabel{rem:4.9}
In \cite[Thm.~2.4]{NO15} we have seen that a $\beta$-periodic function 
$f \: \R \to \R$ is reflection positive if and only if it is of the form 
\[ f(t) = \int_0^\infty e^{-t\lambda} + e^{-(\beta - t)\lambda}\, d\mu_+(\lambda)
\quad \mbox{ for }\quad 0 \leq t \leq \beta, \] 
where $\mu_+$ is a positive measure on $[0,\infty)$. 
For the measure $d\mu(\lambda) = d\mu_+(\lambda) + e^{\beta\lambda} d\mu_+(-\lambda)$ 
on $\R$, this leads to $f = \cL(\mu)$ on $[0,\beta]$. 
The function $f\res_{[0,\beta]}$ is convex and symmetric with respect to 
$\beta/2$, where it has a global minimum. Therefore 
Theorem~\ref{thm:5.9} would only apply to the restriction 
of $f$ to the interval $[-\beta/2,\beta/2]$. 
\end{rem}

For the applications in \cite{JNO16}, we also note the following 
description of reflection negative functions on $\R$: 

\begin{theorem} \mlabel{thm:lk-bernstein} {\rm(Reflection negative functions on 
$(\R,\R_+,-\id_\R)$)}
A symmetric continuous function $\psi \colon \R \to [0,\infty)$ 
is reflection negative with 
respect to $(\R,\R_+,-\id_\R)$ if and only 
if $\psi\res_{(0,\infty)}$ is a Bernstein function. In particular, this is equivalent to 
the existence of $a,b \geq 0$ and a positive measure $\mu$ on $(0,\infty)$ 
with $\int_0^\infty (1 \wedge \lambda)\, d\mu(\lambda) < \infty$ 
such that we have the L\'evy--Khintchine representation 
\[ \psi(t) = a + b|t| + \int_0^\infty (1 - e^{-\lambda |t|})\, d\mu(\lambda). \] 
Here $a,b$ and $\mu$ are uniquely determined by~$\psi$.
\end{theorem}

\begin{prf}  According to \cite[Cor.~3.3]{NO14}, 
a symmetric function $\vphi \: \R \to \R$ is reflection positive 
if and only if $\vphi\res_{\R_+}$ is a Laplace transform of a positive 
measure on $[0,\infty)$ (Theorem~\ref{thm:bern1}). 
Combining this with Schoenberg's Theorem (Remark~\ref{rem:schoen}) and 
Theorem~\ref{thm:lk-bernstein2} proves our assertion.
\end{prf}

\begin{example} \label{ex:3.6} 
For $\alpha \geq 0$, the function 
$\psi(t) := |t|^\alpha$ is reflection negative on $\R$ if and only if 
$0 \leq \alpha \leq 1$ (Examples~\ref{ex:bernstein}(a)). 
\end{example}

\end{document}